\documentclass{article}
\usepackage{longtable}
\usepackage{amsmath}
\usepackage{url}
\usepackage[margin=10truemm]{geometry}
\title{\huge How did Euler solve the equation $xyz(x+y+z) = a$?}
\author{Seiji Tomita and  Oliver Couto}
\date{}
\begin{document}
\maketitle
\begin{abstract}
In this paper, we derived the parametric solution of Euler and Elkies in an elementary manner.\\
In addition we proved there are infinitely many parametric solutions of Euler's and Elkies's family of solutions.

\end{abstract}
\vskip5\baselineskip

\centerline{\huge1. Introduction}
\vskip\baselineskip
According to Elkies\cite{a}, in 1749 Euler takes a look at $xyz(x + y + z) = a$ and says that he has found, with quite some effort.\\
His parametric solution is as follows.
\begin{align*}
x &= \frac{6ast^3(at^4-2s^4)^2}{(4at^4+s^4)(2a^2t^8+10as^4t^4-s^8)}\\[5pt]
y &= \frac{3}{2}\frac{s^5(4at^4+s^4)^2}{t(at^4-2s^4)(2a^2t^8+10as^4t^4-s^8)}\\[5pt]
z &= \frac{2}{3}\frac{(2a^2t^8+10as^4t^4-s^8)}{s^3t(4at^4+s^4)}
\end{align*}

We don't know how he got his solution.

In 2014,  Elkies looked for a solution, and he got some solutions using  algebraic geometry.\\
Simpler solution is  as follows.
\begin{align*}
x &= \frac{(s^4-4a)^2}{2s^3(s^4+12a)}\\[5pt]
y &= \frac{2a(3s^4+4a)^2}{s^3(s4-4a)(s^4+12a)}\\[5pt]
z &= \frac{s(s^4+12a)}{2(3s^4+4a)})
\end{align*}

In 2022,inspired by the problem posted on MathStackExchange\cite{b}, it is related to our problem, 
so we decided to work on elementary derivation of Euler's solution.\\
We proved that there are infinitely many parametric solutions of Euler's family of solutions.\\
In addition we derived Elkies's solution in an elementary manner and proved that there are infinitely many parametric solutions of Elkies's family of solutions.\\
Moreover we showed the small positive solutions table for $a <100$.\\
Furthermore, parametric solution of $wxyz(w+x+y+z)\ =\ a$ was shown in Appendix.\\
An infinitely many solutions are generated from infinite order point of elliptic curve using group law.
In this context, we call generated solutions are family of solutions of infinite order point.

\newpage

\centerline{\huge2. Derive Euler's solution}
\vskip\baselineskip
According to Elkies{}, Euler found a parametric solution of $xyz(x+y+z) = a$ below.
\begin{align*}
x &= \frac{6ast^3(at^4-2s^4)^2}{(4at^4+s^4)(2a^2t^8+10as^4t^4-s^8)}\\[5pt]
y &= \frac{3}{2}\frac{s^5(4at^4+s^4)^2}{t(at^4-2s^4)(2a^2t^8+10as^4t^4-s^8)}\\[5pt]
z &= \frac{2}{3}\frac{(2a^2t^8+10as^4t^4-s^8)}{s^3t(4at^4+s^4)}
\end{align*}

\textbf{We prove that  Euler's family of solutions has an infinitely many parametric solutions where $a$ is arbitrary.}
\vskip\baselineskip

Proof.

$$xyz(x+y+z)\ =\ a$$

We define Euler's solution form as follows.\\
\begin{align*}
x &= \frac{c_1as^s_1t^t_1A^2}{BC}\\[5pt]
y &= \frac{c_2s^s_2t^t_2B^2}{AC}\\[5pt]
z &= \frac{c_3s^s_3t^t_3C}{B}
\end{align*}

Taking $(s1,s2,s3)=(1, 5, -3),\ (t1,t2,t3)=(3, -1, -1),\ (c1,c2,c3)=(6, \frac{3}{2}, \frac{2}{3})$, then
\begin{align*}
x\ &=\ \frac{6Ast^3A^2}{BC}\\
y\ &=\ \frac{3}{2}\frac{s^5B^2}{tAC } \\
z\ &=\ \frac{2}{3}\frac{C}{s^3tB}
\end{align*}

Hence
 $$xyz(x+y+z)\ =\ \frac{a(36As^4t^4A^3+9s^8B^3+4C^2A)}{C^2B}$$

RHS of above equation must be a,then
\begin{equation}
36As^4t^4A^3+9s^8B^3+4C^2A=C^2B\nonumber
\end{equation}

Since $C$ must be rational number then discriminant must be square number.  

Let $U=\frac{A}{B}$ we get quartic curve

$$V^2=-16at^4U^4+4at^4U^3-4s^4U+s^4$$

Quartic is birationally equivalent to an elliptic curve.

$$E: Y^2-4s^2YX+8s^2at^4Y = X^3-4s^4X^2+64at^4s^4X-256at^4s^8$$

E has a point $P(X,Y)=(4s^4,\ 16s^6-8s^2at^4)$.
\vskip\baselineskip
Hence we can obtain a point $2Q(U) = \cfrac{at^4-2s^4}{4at^4+s^4}$\ from $2P(X,Y)$ using group law.
\vskip\baselineskip

Then we obtain
 $$(A,B,C)=(at^4-2s^4, 4at^4+s^4, -s^8+10as^4t^4+2a^2t^8)$$

Finally, we obtain an Euler's solution.
\begin{align*}
x &= \frac{6ast^3(at^4-2s^4)^2}{(4at^4+s^4)(2a^2t^8+10as^4t^4-s^8)}\\[5pt]
y &= \frac{3}{2}\frac{s^5(4at^4+s^4)^2}{t(at^4-2s^4)(2a^2t^8+10as^4t^4-s^8)}\\[5pt]
z &= \frac{2}{3}\frac{(2a^2t^8+10as^4t^4-s^8)}{s^3t(4at^4+s^4)}
\end{align*}

Similarly, we can obtain other new solutions using group law with $P(X,Y)=(4s^4, 16s^6-8s^2at^4)$.

Hence Euler's family of solutions has an infinitely many parametric solutions.

For instance, we can obtain a new solution using $3P(X,Y)$ as follows.
\begin{align*}
x &= \frac{18as^5t^3(-2s^4+at^4)^2(-s^8+10at^4s^4+2a^2t^8)^2(a^2t^8-5s^8+14at^4s^4)}{(at^4+s^4)(a^3t^{12}+3s^4a^2t^8+111at^4s^8+s^{12})(a^6t^{24}+6s^4a^5t^{20}-255s^8a^4t^{16}-790a^3t^{12}s^{12}-2253a^2t^8s^{16}-264s^{20}at^4+s^{24})}\\[5pt]
y &= \frac{1}{6}\frac{-(at^4+s^4)^2(a^3t^{12}+3s^4a^2t^8+111at^4s^8+s^{12})^2(a^2t^8-5s^8+14at^4s^4)}{s^3t(-2s^4+at^4)(-s^8+10at^4s^4+2a^2t^8)(a^6t^{24}+6s^4a^5t^{20}-255s^8a^4t^{16}-790a^3t^{12}s^{12}-2253a^2t^8s^{16}-264s^{20}at^4+s^{24})}\\[5pt]
z &= \frac{2s(a^6t^{24}+6s^4a^5t^{20}-255s^8a^4t^{16}-790a^3t^{12}s^{12}-2253a^2t^8s^{16}-264s^{20}at^4+s^{24})}{(a^2t^8-5s^8+14at^4s^4)t(at^4+s^4)(a^3t^{12}+3s^4a^2t^8+111at^4s^8+s^{12})}\\
\end{align*}
The proof is completed.
\newpage
\centerline{\huge3. Derive Elkies's solution}
\vskip\baselineskip
According to Elkies{}, Elkies found a parametric solution of $xyz(x+y+z) = a$ below.\\
\begin{align*}
x &= \frac{1}{2}\frac{(s^4-4a)^2}{s^3(s^4+12a)}\\[5pt]
y &= \frac{2a(3s^4+4a)^2}{s^3(s^4-4a)(s^4+12a)}\\[5pt]
z &= \frac{1}{2}\frac{s(s^4+12a)}{3s^4+4a}
\end{align*}

\textbf{We show Elkies's family of solutions has an infinitely many parametric solutions where $a$ is arbitrary.}\\
\vskip\baselineskip

Proof.

$$xyz(x+y+z)\ =\ a$$

Taking 
\begin{align*}
x\ &=\frac{1}{2}\frac{A^2}{s^3C}\\[5pt]
y\ &=\frac{2aB^2}{s^3AC}\\[5pt]
z\ &=\frac{1}{2}\frac{sC}{B}
\end{align*}

Hence
 $$xyz(x+y+z)\ =\ \frac{1}{4}\frac{a(A^3B+4aB^3+s^4C^2A)}{s^8C^2}$$

RHS of above equation must be $a$,then\\
\begin{equation}
A^3B+4aB^3+s^4C^2A-4s^8C^2=0\tag{1}
\end{equation}

Since $C$ must be rational number then discriminant must be square number.\\

$$v^2 = -4(A-4s^4)s^4B(A^3+4aB^2)$$

Let $B = -A+4s^4$ then\\ 
\begin{align}
V^2 &= A^3+4aB^2 \nonumber \\
    &= A^3+4aA^2-32s^4aA+64s^8a\tag{2}
\end{align}

In order to find a parametrization for $A, V$, substitute $(A,V)=(s^4+p, s^6+qs^4+rs^2)$ to (2).\\
We obtain $(p,q,r)=(-4a, 0, 12a)$ then $(A,B,C)=(s^4-4a, 3s^4+4a, s^4+12a)$.\\ 
Finally, we obtain Elkies's solution.\\
\begin{align*}
x &= \frac{1}{2}\frac{(s^4-4a)^2}{s^3(s^4+12a)}\\[5pt]
y &= \frac{2a(3s^4+4a)^2}{s^3(s^4-4a)(s^4+12a)}\\[5pt]
z &= \frac{1}{2}\frac{s(s^4+12a)}{3s^4+4a}
\end{align*}

From (2), we define elliptic curve E\\
$$E: V^2 = A^3+4aA^2-32s^4aA+64as^8$$

Since E has a point $P(A,V)=(s^4-4a, s^6+12as^2)$,\\
\vskip\baselineskip
we can obtain a point $2P(A) = \cfrac{1}{4}\cfrac{s^{16}-464s^{12}a+1632s^8a^2+768a^3s^4+256a^4}{s^4(s^4+12a)^2}$.\\
According to Nagell-Lutz theorem, the point $P(A,V)$ is not a point of finite order.\\
Hence Elkies's family of solutions has an infinitely many parametric solutions.\\

We can obtain new solutions using group law with $P(A,V)=(s^4-4a, s^6+12as^2)$.\\
For instance, we can obtain a new solution using $2P(A,V)$ as follows.\\

We obtain\\
\begin{align*}
A &= \frac{1}{4}\frac{(-4a+s^4)(s^{12}-460s^8a-208a^2s^4-64a^3)}{s^4(s^4+12a)^2}\\[5pt]
B &= \frac{1}{4}\frac{(-4a+5s^4)(4a+3s^4)(s^8+56s^4a+16a^2)}{s^4(s^4+12a)^2}\\[5pt]
C &= \frac{-1}{8}\frac{(-16a^2-32s^4a+s^8)(s^{16}+1136s^{12}a-928s^8a^2+1792a^3s^4+256a^4)}{s^8(s^4+12a)^3}
\end{align*}

\begin{align*}
x &= \frac{-1}{4}\frac{(-4a+s^4)^2(s^{12}-460s^8a-208a^2s^4-64a^3)^2}{s^3(s^4+12a)(-16a^2-32s^4a+s^8)(s^{16}+1136s^{12}a-928s^8a^2+1792a^3s^4+256a^4)}\\[5pt]
y &= \frac{-4(s^4+12a)s(-4a+5s^4)^2(4a+3s^4)^2(s^8+56s^4a+16a^2)^2a}{(s^{16}+1136s^{12}a-928s^8a^2+1792a^3s^4+256a^4)(-16a^2-32s^4a+s^8)(-4a+s^4)(s^{12}-460s^8a-208a^2s^4-64a^3)}\\[5pt]
z &= \frac{-1}{4}\frac{(-16a^2-32s^4a+s^8)(s^{16}+1136s^{12}a-928s^8a^2+1792a^3s^4+256a^4)}{s^3(-4a+5s^4)(4a+3s^4)(s^8+56s^4a+16a^2)(s^4+12a)}
\end{align*}
The proof is completed.

\newpage
\centerline{\huge4. Solution table}
\vskip\baselineskip
\begin{center}
Small positive solutions by brute force search with $a<100$.\\
\end{center}

\begin{longtable}[c]{rrrr}
\caption{Solutions of $xyz(x+y+z) = a$}
\label{longtablesample} \\
\hline
a & x & y & z \\
\hline \endhead
 1 & 3/2 & 4/3 & 1/6 \\
 2 & 5/2 & 5/6 & 4/15 \\
 3 & 1 & 1 & 1 \\
 4 & 7/2 & 36/35 & 7/30 \\
 5 & 4 & 1/2 & 1/2 \\
 6 & 2 & 3/2 & 1/2 \\
 7 & 10/3 & 21/20 & 5/12 \\
 8 & 2 & 1 & 1 \\
 9 & 2 & 2 & 1/2 \\
 10 & 5/2 & 4/3 & 2/3 \\
 11 & 2 & 11/6 & 2/3 \\
 12 & 5/3 & 27/20 & 5/4 \\
 13 & 3/2 & 3/2 & 4/3 \\
 14 & 2 & 4/3 & 7/6 \\
 15 & 3 & 1 & 1 \\
 16 & 3 & 8/3 & 1/3 \\
 17 & 5/2 & 10/7 & 34/35 \\
 18 & 4 & 3/2 & 1/2 \\
 19 & 3 & 3 & 1/3 \\
 20 & 2 & 2 & 1 \\
 21 & 5/2 & 5/2 & 3/5 \\
 22 & 9/2 & 25/6 & 2/15 \\
 23 & 5 & 9/10 & 23/30 \\
 24 & 4 & 1 & 1 \\
 25 & 4 & 9/4 & 5/12 \\
 26 & 4 & 2 & 1/2 \\
 27 & 9/4 & 25/12 & 16/15 \\
 28 & 3 & 7/3 & 2/3 \\
 29 & 4 & 29/20 & 4/5 \\
 30 & 4 & 15/4 & 1/4 \\
 31 & 8/3 & 31/12 & 3/4 \\
 32 & 5 & 5/3 & 8/15 \\
 33 & 2 & 2 & 3/2 \\
 34 & 7/2 & 7/3 & 9/14 \\
 35 & 3 & 3/2 & 4/3 \\
 36 & 3 & 2 & 1 \\
 37 & 20 & 37/15 & 1/30 \\
 38 & 6 & 12/7 & 19/42 \\
 39 & 4 & 3/2 & 1 \\
 40 & 3 & 3 & 2/3 \\
 41 & 3 & 25/12 & 16/15 \\
 42 & 15/4 & 7/5 & 5/4 \\
 43 & 9/2 & 2 & 2/3 \\
 44 & 7 & 7/4 & 11/28 \\
 45 & 4 & 3 & 1/2 \\
 46 & 8 & 25/12 & 4/15 \\
 47 & 25/6 & 32/15 & 3/4 \\
 48 & 2 & 2 & 2 \\
 49 & 6 & 3/2 & 2/3 \\
 50 & 5 & 5/2 & 1/2 \\
 51 & 4 & 17/12 & 4/3 \\
 52 & 6 & 26/15 & 3/5 \\
 53 & 7/2 & 7/2 & 4/7 \\
 54 & 5 & 8/5 & 9/10 \\
 55 & 5 & 11/3 & 1/3 \\
 56 & 4 & 2 & 1 \\
 57 & 16/5 & 25/8 & 4/5 \\
 58 & 20/3 & 6/5 & 5/6 \\
 59 & 5 & 49/30 & 20/21 \\
 60 & 8 & 3/2 & 1/2 \\
 61 & 32/3 & 9/2 & 1/12 \\
 62 & 9/4 & 25/12 & 31/15 \\
 63 & 3 & 3 & 1 \\
 64 & 9 & 1 & 2/3 \\
 65 & 3 & 13/6 & 3/2 \\
 66 & 5/2 & 5/2 & 8/5 \\
 67 & 15/2 & 32/15 & 5/12 \\
 68 & 4 & 4 & 1/2 \\
 69 & 13/4 & 13/4 & 23/26 \\
 70 & 6 & 35/6 & 1/6 \\
 71 & 10 & 49/20 & 8/35 \\
 72 & 16 & 25/8 & 3/40 \\
 73 & 6 & 6 & 1/6 \\
 74 & 15/2 & 10/3 & 4/15 \\
 75 & 4 & 5/2 & 1 \\
 76 & 7/2 & 18/7 & 7/6 \\
 77 & 5 & 16/5 & 11/20 \\
 78 & 4 & 13/4 & 3/4 \\
 79 & 8 & 25/12 & 9/20 \\
 80 & 5 & 2 & 1 \\
 81 & 9/2 & 4 & 1/2 \\
 82 & 10/3 & 5/2 & 41/30 \\
 83 & 5 & 5/2 & 4/5 \\
 84 & 3 & 2 & 2 \\
 85 & 9 & 16/9 & 17/36 \\
 86 & 3 & 8/3 & 3/2 \\
 87 & 9 & 25/3 & 1/15 \\
 88 & 3 & 3 & 4/3 \\
 89 & 27/2 & 2/3 & 2/3 \\
 90 & 4 & 2 & 3/2 \\
 91 & 4 & 9/4 & 4/3 \\
 92 & 5 & 23/5 & 2/5 \\
 93 & 11/2 & 11/2 & 3/11 \\
 94 & 9/2 & 2 & 4/3 \\
 95 & 5 & 4 & 1/2 \\
 96 & 4 & 3 & 1 \\
 97 & 10/3 & 49/20 & 45/28 \\
 98 & 7/2 & 7/2 & 1 \\
 99 & 5/2 & 5/2 & 11/5 \\

\end{longtable}

\label{longtable}
\newpage
\centerline{\huge5. Final remarks}
\vskip\baselineskip
We used the Euler's solution form, combination of $(A,B,C)$ to derive Euler's parametric solution.\\
So, we think there must be parametric solutions other than Euler's form.\\
The same could be said of Elkies case in a similar way.
Furthermore, it might be interesting to construct the parametric solutions using $(A,B,C,D)$.

\vskip5\baselineskip
\centerline{\huge Appendix}
\vskip\baselineskip
We consider the extension of $xyz(x+y+z) = a$ and show $wxyz(w+x+y+z) = a$ has infinitely many  parametric solutions where $a$ is arbitrary.
The problem $wxyz(w+x+y+z) = 1$ was posted on mathoverflow.net\cite{c}.
\vskip\baselineskip
$\textbf{ An equation wxyz(w+x+y+z) = a}$ \textbf{has an infinitely many  parametric solutions where $a$ is arbitrary.}\\

Proof.

$$wxyz(w+x+y+z)=a$$

\begin{align*}
w &= \frac{ac_1t^{t_1}}{ABC}\\
x &= \frac{c_2Bt^{t_2}}{A}\\
y &= \frac{c_3At^{t_3}}{C}\\
z &= \frac{c_4Ct^{t_4}}{B}
\end{align*}

Hence
 $$wxyz(w+x+y+z) = \frac{c_1t^{t_1}c_2t^{t_2}c_3t^{t_3}c_4t^{t_4}(ac_1t^{t_1}+c_2B^2t^{t_2}C+c_3A^2t^{t_3}B+c_4C^2t^{t_4}A)}{A^2B^2C^2}$$ 

RHS of above equation must be a,then
\begin{align*}
&(c_1c_2c_3c_4^2t^{t_1+t_2+t_3+2t_4}A-A^2B^2)C^2\\
&+c_1c_2^2c_3c_4t^{t_4+t_3+t_1+2t_2}B^2C\\
&+ac_1^2c_2c_3c_4t^{t_4+t_3+t_2+2t_1}+c_1c_2c_3^2c_4t^{t_1+t_2+t_4+2t_3}A^2B=0
\end{align*}

Since $C$ must be rational number then discriminant must be square number.  
\begin{align*}
v^2 &= 4c_1c_2c_3^2c_4t^{t_4+t_2+t_1+2t_3}B^3A^4\\
    &-4c_1^2c_2^2c_3^3c_4^3Bt^{3t_4+3t_3+2t_2+2t_1}A^3\\
    &+a4c_1^2c_2c_3c_4t^{t_4+t_3+t_2+2t_1}B^2A^2\\
    &-a4c_1^3c_2^2c_3^2c_4^3t^{3t_4+2t_3+2t_2+3t_1}A\\
    &+c_1^2c_2^4c_3^2c_4^2B^4t^{2t_4+2t_3+4t_2+2t_1}
\end{align*}

Obviously, quartic is birationally equivalent to an elliptic curve.  

For instance, take $(c_1,c_2,c_3)=(1,1,1),(t_1,t_2,t_3)=(1,1,1),B=1$, then  

$$v^2 = 4t^5A^4-4t^{10}A^3+4t^5aA^2-4t^{10}aA+t^{10}$$

Quartic is transformed to an elliptic curve below.

$$E: Y^2-4t^5aYX-8t^{15}Y = X^3+(4t^5a-4t^{10}a^2)X^2-16t^{15}X-64t^{20}a+64t^{25}a^2$$

$E$ has a point $P(X,Y)=(-4t^5a+4t^{10}a^2,\ -16t^{10}a^2+16t^{15}a^3+8t^{15})$.\\  
We obtain $$2P(X)= \frac{4((-a^2-2a^5+a^8)t^{20}+(-2a^4-a-4a^7)t^{15}+(1+6a^6+6a^3)t^{10}+(-2a^2-4a^5)t^5+a^4)}{((2a^3+1)t^5-2a^2)^2}$$  
According to Nagell-Lutz theorem, the point $P(X,Y)$ is not a point of finite order, hence we can obtain infinitely many parametric solutions.  

Thus we can obtain a quartic point $$2Q(A) = \frac{-t^5(-2a^2+2t^5a^3+t^5)}{t^{10}a^4-2t^5a^3-t^5+a^2}$$  

Then we obtain $(A,C)$  
\begin{align*}
A &= \frac{-t^5(-2a^2+2t^5a^3+t^5)}{t^{10}a^4-2t^5a^3-t^5+a^2}\\[5pt]
C &= \frac{t^{10}+a-2t^5a^2+t^{10}a^3}{t^{10}a^2-1}
\end{align*}

Finally we obtain $(w,x,y,z)$
\begin{align*}
w &= \frac{-a(t^{10}a^4-2t^5a^3-t^5+a^2)(t^5a+1)(t^5a-1)}{t^4(-2a^2+2t^5a^3+t^5)(t^{10}+a-2t^5a^2+t^{10}a^3)}\\[5pt]
x &= \frac{-(t^{10}a^4-2t^5a^3-t^5+a^2)}{t^4(-2a^2+2t^5a^3+t^5)}\\[5pt]
y &= \frac{-t^6(-2a^2+2t^5a^3+t^5)(t^5a+1)(t^5a-1)}{(t^{10}a^4-2t^5a^3-t^5+a^2)(t^{10}+a-2t^5a^2+t^{10}a^3)}\\[5pt]
z &= \frac{(t^{10}+a-2t^5a^2+t^{10}a^3)t}{(t^5a+1)(t^5a-1)}
\end{align*}

In this way we can obtain infinitely many parametric solutions.\\
The proof is completed.

\vskip\baselineskip

\end{document}